\documentclass[12pt]{amsart}
\usepackage{epsfig}
\newcommand {\supplus}{\mathop{{\supset}\llap{\raise 
0.5pt\hbox{\normalfont\small+}\hskip 0.5pt}}} 
\newcommand {\subplus}{\mathop{{\subset}\llap{\raise 
0.5pt\hbox{\normalfont\small+}\hskip 0.5pt}}}  
\newcommand {\Cee}    {{\mathbb  C}}

\newcommand {\Kee}    {{\mathbb  K}}

\newcommand {\Zee}    {{\mathbb  Z}}

\newcommand {\faut}   {{\mathfrak{aut}}} 
\newcommand {\fb}     {{\mathfrak{b}}}

\newcommand {\fder}   {{\mathfrak{der}}}   %

\newcommand {\fd}     {{\mathfrak{d}}}

\newcommand {\fg}     {{\mathfrak{g}}}    %
\newcommand {\fgl}    {{\mathfrak{gl}}}  %
\newcommand {\fh}     {{\mathfrak{h}}}

\newcommand {\fo}     {{\mathfrak{o}}}
\newcommand {\fosp}   {{\mathfrak{osp}}}
\newcommand {\fp}    {{\mathfrak{p}}}   %
\newcommand {\fpe}    {{\mathfrak{pe}}}   %

\newcommand {\fpsl}   {{\mathfrak{psl}}}

\newcommand {\fs}     {{\mathfrak{s}}}

\newcommand {\fsh}    {{\mathfrak{sh}}}

\newcommand {\fsl}    {{\mathfrak{sl}}}

\newcommand {\fsp}    {{\mathfrak{sp}}}
\newcommand {\fspe}   {{\mathfrak{spe}}}

\newcommand {\fsvect} {{\mathfrak{svect}}}

\newcommand {\fvect}  {{\mathfrak{vect}}}   %

\newcommand {\cal} {\mathcal}

\newcommand {\cI}     {{\cal I}}

\newcommand {\cL}     {{\cal L}}

\def \opname#1#2%
  {\expandafter\newcommand \csname #1\endcsname {{\mathop{#2}\nolimits}}}

\newcommand{\rmname}[1]
  {\expandafter\newcommand \csname #1\endcsname {{\operatorname{#1}}}}

\newcommand{\rmnameii}[2]
  {\expandafter\newcommand \csname #1\endcsname {{\operatorname{#2}}}}

\rmname{act}
\rmname{Ad}
\rmname{Add}
\rmname{ad}
\rmname{Alt}
\rmname{alt}
\rmname{Ann}
\rmname{antidiag}
\rmname{Ber}
\rmname{ber}
\rmname{Br}
\rmname{card}
\rmname{ch}
\rmname{Char}
\rmname{cem}
\rmname{cj}
\rmname{Cliff}
\rmname{cntr}
\rmname{codim}
\rmname{coind}
\rmname{const}
\rmname{col}
\rmname{cork}
\rmname{cpr}
\rmname{diag}
\rmnameii{Div}{div}
\rmname{Def}
\rmname{Der}
\rmname{Dim}
\rmname{End}
\rmname{Even}
\rmname{Ext}
\rmname{gr}
\rmname{Hom}
\rmname{HT}
\rmnameii{Ht}{ht}
\rmname{hwt}
\rmname{Id}
\rmname{id}
\rmname{ind}
\rmname{Ind}
\rmname{Inf}
\rmname{irr}
\rmname{Le}
\rmname{Lie}
\rmname{lwt} 
\rmname{mult}
\rmname{Mor}
\rmname{nm}
\rmname{Ob}
\rmname{Odd}
\rmname{Osc}
\rmname{per}
\rmname{Pic}
\rmname{pr}
\rmname{pro}
\rmname{Prime}
\rmname{Proj}
\rmname{prt}
\rmname{pt}
\rmname{Q}
\rmname{qet}
\rmname{qtr}
\rmname{rd}
\rmname{rk}
\rmname{row}
\rmname{Res}
\rmname{salt}
\rmname{Sch}
\rmname{SBr}
\rmname{scalar}
\rmname{Ser}
\rmname{sign}
\rmname{Smbl}
\rmname{spin}
\rmname{ssym}
\rmname{str}
\rmname{st}
\rmname{sgn}
\rmname{sq}
\rmname{symm}
\rmname{supp}
\rmname{Supp}
\rmname{St}
\rmname{Spec}
\rmname{Spm}
\rmname{tr}
\rmname{vpt}
\rmname{weyl}
\rmname{Weyl}
\rmname{Witt}

\opname{vvol}  {{v\hspace{-0.1ex}o\hspace{-0.02ex}l\/}}
\opname{pnt}  {\text{\normalfont pt}}
\opname{Span} {{Span}}
\opname{slim} {\overline{\lim}}
\opname{Vol}  {{V\hspace{-0.55ex}o\hspace{-0.02ex}l\/}}
\opname{QVol} {{Q\hspace{-0.3ex}V\hspace{-0.55ex}o\hspace{-0.02ex}l\/}}
\opname{PoVol}{{P\hspace{-0.35ex}o\hspace{-0.25ex}V\hspace{-0.55ex}o\hspace{-0.02ex}l\/}}
\opname{BVol} {{B\hspace{-0.2ex}V\hspace{-0.55ex}o\hspace{-0.02ex}l\/}}
\opname{Par}  {{P\hspace{-0.3ex}a\hspace{-0.05ex}r\/}}

\rmname{Mat}
\rmname{Bil}
\rmname{Diff}
\rmname{Ker}
\rmname{Herm}
\rmname{Coker}
\rmname{Conn}
\rmname{Covect}
\rmname{Vect}
\rmname{Int}

\rmnameii {IM} {Im}
\rmnameii {RE} {Re}

\opname{Aut} {{A\hspace{-0.2ex}u\hspace{-0.1ex}t\/}}
\opname{GL} {{G\hspace{-0.3ex}L}}
\opname{SL} {{S\hspace{-0.3ex}L}}
\opname{Exp} {{E\hspace{-0.2ex}x\hspace{-0.1ex}p\/}}
\opname{GQ} {{G\hspace{-0.2ex}Q}}
\opname{OSp} {{O\hspace{-0.25ex}S\hspace{-0.15ex}p\/}}
\opname{Out} {{O\hspace{-0.25ex}u\hspace{-0.15ex}t\/}}
\opname{Spp} {{S\hspace{-0.2ex}p\/}}
\opname{SpO} {{S\hspace{-0.2ex}p\hspace{-0.02ex}O\/}}
\opname{Pe} {{P\hspace{-0.25ex}e\/}}
\opname{SPe} {{S\hspace{-0.25ex}P\hspace{-0.25ex}e\/}}
\opname{Spin} {{S\hspace{-0.25ex}p\hspace{-0.05ex}i\hspace{-0.1ex}n\/}}
\opname{Iso} {{I\hspace{-0.25ex}s\hspace{-0.1ex}o\/}}
\opname{SSPe} {{S\hspace{-0.25ex}S\hspace{-0.15ex}P\hspace{-0.25ex}e\/}}
\opname{PeU} {{P\hspace{-0.25ex}e\hspace{-0.1ex}U\/}}
\opname{QU} {{Q\hspace{-0.15ex}U\/}}
\opname{U} {{U\/}}

\opname{cGQ} {{\cal G \hspace{-0.2em} Q \/}}
\opname{cSL} {{\cal S \hspace{-0.2em} L \/}}
\opname{cGL} {{\cal G \hspace{-0.2em} L \/}}
\opname{cOSp} {{\cal O \hspace{-0.2em} S \hspace{-0.3em} \it p\/}}
\opname{cPe} {{\cal P \hspace{-1.5pt} \it e\/}}
\opname{cVect} {{\cal V \hspace{-1.5pt} \it e\hspace{-0.1ex}c\hspace{-0.1ex}t\/}}
\opname{cVol} {{\cal V \hspace{-1.5pt} \it o\hspace{-0.1ex}l\/}}
\opname{cAut} {{\cal A \hspace{-0.2em} \it u\hspace{-0.1em}t\/}}
\opname{cCovect} {{\cal C \hspace{-1.5pt}
     \it o\hspace{-0.1ex}v\hspace{-0.1ex}e\hspace{-0.1ex}c\hspace{-0.1ex}t\/}}
\opname{CW} {{C\hspace{-0.15ex}W}}

\opname {Ab}   {{\sf Ab}}
\opname {Alg}   {{\sf Alg}}
\opname {ASch}  {{\sf Aff\;Sch}}
\opname {Funct}   {{\sf Funct}}
\opname {Gr}   {{\sf Gr}}
\opname {Grf}  {{{\sf Gr}_f}}
\opname {Mods}   {{\sf mods}}
\opname {Rings}   {{\sf Rings}}
\opname {Salg}   {{\sf Salg}}
\opname {Sets} {{\sf Sets}}
\opname {SSMan} {{\sf SMan}}
\opname {Top}  {{\sf Top}}
\opname {Vebun}   {{\sf Vebun}}

\newcommand {\ev} {{\bar0}}
\newcommand {\od} {{\bar1}}
\newcommand {\eps} {\varepsilon}

\newcommand {\tto} {\longrightarrow}
\newcommand {\pder}[1] {{\frac{\partial}{\partial {#1}}}}
\newcommand {\pderf}[2] {{\frac{\partial {#1}}{\partial {#2}}}}

\newcommand {\bcdot}   {\mathbin{\hbox{\raise.4ex\hbox{\bf.}}}} % bold \cdot

\newcommand {\secno} {}
\newcommand {\ssecfont} {\normalfont\bf}

\newtheorem{Theorem}{\secno Theorem}

\newtheorem{Lemma}[Theorem]{\secno Lemma}

\newenvironment {th*}[1]
    {\gdef\thname{#1} \begin{thn}}%
    {\end{thn}}
\newtheorem{thn}[Theorem] {\thname}

\theoremstyle{definition}

\newenvironment {ex*}[1]
    {\gdef\thname{#1} \begin{exn}}%
    {\end{exn}}
\newtheorem{exn}[Theorem]{\thname}

\theoremstyle{remark}

\newenvironment {rem*}[1]
    {\gdef\thname{#1} \begin{remn}}%
    {\end{remn}}
\newtheorem{remn}[Theorem]{\thname}

\newcommand {\ssec}{\subsection*}

\newcommand {\ssbegin}[2]
  {\def \secno {\gdef \secno {}{\ssecfont #1. }}%
   \begin{#2}}
\begin{document}

\title{Indecomposable representations of Lie superalgebras}

\author{Dimitry Leites}

\thanks{Financial support of NFR and technical of V.~Serganova, 
A.~Shapovalov and P.~Grozman, as well as a crucial hint of A.~A.~Kirillov, 
are gratefully acknowledged.}



\begin{abstract} 
In 1960's I.~Gelfand posed a problem: describe indecomposable 
representations of any simple infinite dimensional Lie algebra $\fg$ 
of polynomial vector fields.  Here by applying the elementary 
technique of Gelfand and Ponomarev a toy model of the problem is 
solved: finite dimensional indecomposable representations of $\fvect 
(0|2)$, the Lie superalgebra of vector fields on the 
$(0|2)$-dimensional superspace, are described.

Since $\fvect (0|2)$ is isomorphic to $\fsl(1|2)$ and $\fosp(2|2)$, 
their representations are also are described.  The result is 
generalized in two directions: for $\fsl(1|n)$ and $\fosp(2|2n)$.  
Independently and differently J.~Germoni described indecomposable 
representation of the series $\fsl(1|n)$ and several individual Lie 
superalgebras.

Partial results for other simple Lie superalgebras without Cartan 
matrix are reviewed.  In particular, it is only for 
$\fvect(0|2)\simeq\fsl(1|n)$ and $\fs\fh(0|4)\simeq\fpsl(2|2)$ that 
the typical irreducible representations can not participate in 
indecomposable modules; for other simple Lie superalgebras without 
Cartan matrix (of series $\fvect(0|n)$, $\fsvect(0|n)$, 
$\widetilde\fsvect(0|n)$, $\fs\fpe(n)$ for $n\geq 3$ and $\fs\fh(0|m)$ 
for $m\geq 5$) one can construct indecomposable representations with 
arbitrary composition factors.

Several tame open problems are listed, among them a description of odd 
parameters, previously ignored.
\end{abstract}

\address{Department of Mathematics, University
of Stockholm, Roslagsv.  101, Kr\"aftriket hus 6, SE-106 91,
Stockholm, Sweden; mleites@mathematik.su.se}

\maketitle

{\it To the memory of Misha Saveliev} together with whom we began to 
apply the description of representations of $\fosp(N|2)$ to the 
solution of $N$-extended Leznov-Saveliev equations, cf.  \cite{LSS}.

\section*{Introduction}

This is a version of the paper published in: In: Sissakian A.~N.
(ed.), {\it Memorial volume dedicated to Misha Saveliev and Igor
Luzenko}, JINR, Dubna, 2000, 126--131.  

Meanwhile Noam Shomron \cite{Sh} rediscovered an unpublished
A.~Shapovalov's results on modules over $\fvect(0|m)$, $\fsvect(0|m)$,
$\fpe(m)$, $\fs\fpe(m)$ and $\fh(0|m)$, $\fsh(0|m)$.  Shomron
formulated it only for $\fvect(0|m)$ but went further in formulation
and his proof is impressive.

In what follows the ground field is $\Cee$.  This paper is an attempt 
to review the subject, in particular, formulate tame open problems.

\ssec{0.0.  Prehistory} The study of irreducible finite dimensional 
representations of simple finite dimensional Lie algebras over $\Cee$ 
is a natural problem.  It turns out that such representations are 
completely reducible and, therefore, it suffices to study irreducible 
modules.  Situation changes when we consider infinite dimensional 
modules, even if the module is ``semi-infinite'' in a sence, e.g., 
possesses a highest or lowest weight vector.  Among such 
representation and their ``infinite in both ways'' generalizations, 
Harish-Chandra modules, there is no complete reducibility but the 
problem of description of all indecomposable modules in these 
categories seem to be wild.  (So we can describe either finite 
dimensional ({\it hence}, irreducible), or infinite dimensional {\it and} 
irreducible, or nothing.)

The study of invariant differential operators on manifolds is closely 
related with continuous in a natural topology infinite dimensional 
representations of simple infinite dimensional Lie algebras of vector 
fields $\cL$.  Though the representations themselves are no longer 
completely reducible, the description of irreducible ones reduces, to 
an extent, to finite dimensional representations of simple finite 
dimensional Lie algebras $L_{0}$, the linear parts of $\cL$ (see 
\cite{R} and an elucidation \cite{BL}).  This reduction to finite 
dimensional representation of simple Lie algebras was, perhaps, a 
motivation for I.~Gelfand to consider the classification problem of 
indecomposable representations (say, in the class of modules with 
highest or lowest weight vectors) of such Lie algebras not hopeless 
(at least, in the particular and exceptional case of $\fvect(1)=\fder 
\,\Cee[x]$).  Observe that this problem is still open.

Passing to Lie superalgebras, we are {\it forced} to consider their 
indecomposable representations even in the above problems for several 
reasons: 

(1) even finite dimensional representations of simple Lie 
superalgebras are never completely reducible (with the only exception 
of $\fosp(1|2n)$; the proof of this fact is similar to that for Lie 
algebras, cf.  \cite{B});

(2) the description of irreducible and continuous in a natural 
topology infinite dimensional representations of simple infinite 
dimensional Lie superalgebras of vector fields $\cL$ does not reduce 
to finite dimensional representations of the corresponding finite 
dimensional Lie superalgebras $L_{0}$ but even if we confine ourselves 
to such representations we have to face problem (1).  (This was the 
reason problem (1) was first tackled; G.~Shmelev \cite{Sh1}, 
\cite{Sh2} considered problem (1) only in conneciton with problem (2) 
in a particular case of the latter: for $L_{0}=\fosp(m|2n)$.)

One more reason is provided by the necessity to classify embeddings 
of $\fosp(N|2)$ into the simple Lie superalgebras; this classification 
is a part of the explicit solution of vector-valued generalizations 
of $N$-extended Leznov-Saveliev equations considered, so far, for 
$N=1$ only, see \cite{LSS}.

\ssec{0.1.  History.  Main result} Unless otherwise stated, in what 
follows $\fg$ is either $\fsl(1|n)$ or $\fosp(2|2n)$.  For $\fg$, the 
description of {\it irreducible} finite dimensional representations is 
complete, cf.  \cite{BL}, \cite{L1}.  The next step --- a description 
of indecomposable modules --- was performed under the assumption of 
$\fh$-diagonalizability in \cite{Sh1}, \cite{Sh2}, \cite{L2} who used 
a key observation A.~A.~Kirillov made: he related the problem with a 
result of \cite{ZN}.

A draft of this paper was written in 1989; I delayed the publication 
because I wanted to compute the odd parameters.  Regrettably, it is 
still an open problem.  Meanwhile the problem considered in this paper 
was solved for $\fsl(1|2)$ by Su Yucai \cite{S2} (with minor 
omissions).  He also classifieded indecomposable infinite dimensional 
(Harish-Chandra) modules over $\fosp(1|2)$, see \cite{S1}.

Here we complete the description of indecomposable finite dimensional 
representations of $\fg$ in full generality.  Our results are based on 
the deep results of I.~Gelfand and Ponomarev \cite{GP}.  Not only are 
they deep, they are obtained by elementary methods and I consider 
this as an advantage.  

Meanwhile J.~Germoni independently considered 
the case $\fsl(1|n)$ by much more sophisticated methods (which enabled 
him to relate the result with quivers, etc.), cf.  \cite{G1}, and 
completely classified indecomposable representation.  However, 
\cite{G1} and its continuation \cite{G2} do not mention either 
$\fosp(2|2n)$ considered in \cite{L2} or odd parameters discussed here.

Observe that the description of indecomposable representations for 
$\fsl(1|n)$ and $\fosp(2|2n)$ are obtained by practically identical 
means, be they either elementary (as here) or more involved (as in 
\cite{G1}).

\ssec{0.2. A related problem}  For $\fgl(m|n)$ or $\fsl(m|n)$, Sergeev 
proved \cite{Se} that any finite dimensional representation realizable 
in the space of tensors $T(V\oplus\Pi(V))$ is completely reducible.  
There is no complete reducibility for representations realizable in 
the general tensor algebra $T(V\oplus\Pi(V))\oplus 
V^{*}\oplus\Pi(V^{*}))$.  To distinguish a subalgebra of the tensor 
algebra inside of which the complete reducibility takes place for 
$\fg\neq \fgl$ is an interesting open problem for various 
representations $V$, the identity representation of the matrix 
algebra, the representations of minimal dimention and the adjoint 
representation are most interesting to consider.

\ssec{0.3.  On tame and wild representations} The problem of 
description of indecomposable representations of $\Lambda (n)$ for $n>2$ is a 
wild one (this inference from \cite{GP} is made in \cite{ZN}).  
Similarly, one deduces that to describe indecomposable representations 
of the direct sum of several copies of $\Lambda (2)$ is a wild problem.

Germoni proved that the description of indecomposable representations 
of $\fg=\fsl(p|q)$, where $1<p\leq q$, is a wild problem \cite{G1} by 
cohomology arguments.  Roughly speaking, it seem to be related with 
the possibility to embed the direct sum of two copies of $\fsl(1|1)$ 
into $\fg$ since this reduces to description of indecomposable 
representations of the direct sum of several copies of $\Lambda (2)$.
The arguments are, however, subtler, and it is an open problem to 
consider the indecomposable representations of simple (and close to 
them) Lie superalgebras of other types.  

For certain particular representation of any superalgebra one can 
always obtain the final result.  For example, P.~Grozman recently 
computed with the help of his package SuperLie, see \cite{GL}, the 
decomposition of $\fg\otimes \fg$ for certain particular Lie 
superalgebras $\fg$.  In \cite{G2} Germoni considered indecomposable 
representations of two of the exceptional Lie superalgebras.

\ssec{0.4.  Miscellanea} Several such partial results 
obtained with P.~Grozman constitute explicit decomposition of 
$\fg\otimes\fg$ for $\fg=\fsl(1|1)$, $\fsl(1|2)$, $\fp\fsl(2|2)$ and 
some other algebras.  They are prompted by A.~Vaintrob who hopes to 
relate them with the description of invariants of links and tangles 
\`a la \cite{BS} and are presented elsewhere.

\section*{\protect \S 1. Representations of $\fsl(1|1)$}

Thanks to the super analog of I.~Schur lemma, to consider the 
irreducible representations of $\fsl(1|1)$ is the same as to consider 
same of the (associative) Clifford superalgebra 
$Cl_\hbar(2)=U(\fsl(1|1))/(E-\hbar)$ for $\hbar\in \Cee$, where $E$ is 
the central element of $\fsl(1|1)$.  General 
Algebra and Linear Algebra in Superspaces yield the following 
statement:

\ssbegin{1.1}{Proposition} The irreducible representations of 
$Cl_\hbar(2)$ are:

{\em 1)} for $\hbar\neq 0$: the $1|1$-dimensional modules $V^\hbar$ 
with even highest vector on which $E$ acts as the operator of 
multiplication by the scalar $\hbar$ and $\Pi(V^\hbar)$;

{\em 2)} for $\hbar=0$: the $1|0$-dimensional trivial module ${\bf 1}$ 
and $\Pi({\bf 1})$.
\end{Proposition}

{\bf Problem} Describe representations of $\U(\fsl(1|1))/(E^n-\hbar)$.

Observe that for $\hbar=0$ the superalgebra $Cl_\hbar(2)$ turns into 
the Grassmann superalgebra $\Lambda (2)$ with two generators.  In order to 
study indecomposable representations of $\Lambda (n)$, recall some general 
results on $\Lambda (n)$-modules.

Modules over Grassmann superalgebras are described in \cite{BG}.  Let 
$\Lambda (n)$ be Grassmann superalgebra with $n$ indeterminates, i.e., the 
free supercommutative superalgebra with $n$ odd indeterminates 
$\theta_1, \dots , \theta_n$.  Setting $\deg \theta_i=1$ we define a 
$\Zee$-gradation on $\Lambda (n)$.  In this section we only consider left 
unitary $\Zee$-graded modules with a {\it compatible} 
$\Zee/2$-grading, i.e., the parity of the elements is their degree 
modulo 2 and $\Lambda (n)_i\otimes V_j\subset V_{i+j}$ for $i, j\in\Zee $.  
Given a $\Lambda (n)$-module $V=\oplus V_i$ and $r\in\Zee $, define $V[r]$ by 
setting $V[r]_i=V_{r-i}$.

A $\Lambda (n)$-module $V$ is called {\it reduced} if $\Theta V=0$ for 
$\Theta=\theta_1\cdot\dots \cdot\theta_n$ (the highest term in 
$\Lambda (n)$).

\begin{Lemma}{\em (\cite{BG})} {\em i)} Any free $\Lambda (n)$-module 
$V$ is of the form $F=\oplus_r \Lambda (n)[r]$.

{\em ii)} Any $\Lambda (n)$-module $V$ can be represented in the form
$F\oplus V^{rd}$, where $F$ is free and $V^{rd}$ is reduced.

{\em iii)} Any indecomposable $\Lambda (n)$-module is either isomorphic to 
$\Lambda (n)[r]$ for some $r$ or is reduced.
\end{Lemma}

Since we are only interested in compatibility of the action with 
parity, it suffices to consider shifts of grading modulo 2 in which 
case they coinside with the change of parity functor $\Pi$.

\ssec{1.3.  Reduced modules over $\Lambda (2)$} Let $a$ and $b$ be generators 
of $\Lambda (2)$.  On a reduced $\Lambda (2)$-module $V$ the operators $a$ and $b$ 
satisfy
$$
a^2=b^2=ab=ba=0
$$
and, therefore, results from \cite{GP} are applicable.  Let us recall 
these results (and adjust them to supercase).  We will depict 
$\Lambda (2)$-modules by directed graphs whose nodes stand for subsuperspaces 
of $V$; the action of $a$ is depicted by a horisontal arrow, that of 
$b$ by a vertical one; up to changes of parity both actions are 
isomorphisms representable in a basis by the identity matrix.  The 
curvy arrow may stand for either $a$ or $b$ and corresponds to the 
Jordan cell with eigenvalue $\mu\in\Cee $.

{\bf Type I modules} $V(p+q\varepsilon; dir)$: are determined by their 
dimension $p+q\varepsilon$, where $q=p$ or $q=p\pm 1$, and the 
direction ($dir=in$ or $dir=out$) indicating the submodule (each node 
is $1$-dimensional).  Here is the diagram representing 
$V(p+q\varepsilon; dir)$ for $(p, q)=(3, 2)$ and $(2, 3)$, 
respectively:
$$
\begin{matrix}
    \circ&\tto&\circ&&\cr
    &&\uparrow&&\cr
    &&\circ&\tto&\circ\cr
    &&&&\uparrow\cr
    &&&&\circ\cr
    \end{matrix} \quad \quad \text{and}\quad \quad \begin{matrix}
    \circ&&&&\cr
    \uparrow&&&&\cr
    \circ&\tto&\circ&&\cr
    &&\uparrow&&\cr
    &&\circ&\tto&\circ\cr
    \end{matrix} 
$$

{\bf Type II modules} $V(p; m+n\varepsilon; dir; \mu)$ are 
generalizations of Type I modules by means of the operator depicted by 
a curvy arrow (here depicted by the composition of arrows of two solid 
lines) from the space represented by the most left and upper node to 
the most right and low one or the other way round.
$$
\begin{matrix}
    \circ&\tto&\circ&&&&\cr
    \Downarrow&&\uparrow&&&&\cr
    \Downarrow&&\circ&\tto&\circ&&\cr
    \Downarrow&&&&\uparrow&&\cr
    \Downarrow&&&&\circ&\tto&\circ\cr
    &\Rightarrow&\Rightarrow&\Rightarrow&\Rightarrow&\Rightarrow
    &\Uparrow\end{matrix} \quad \quad \text{and}\quad \quad \begin{matrix}
    \circ&\Leftarrow&\Leftarrow&\Leftarrow&\Leftarrow&\cr
    \uparrow&&&&&\Uparrow\cr
    \circ&\tto&\circ&&&\Uparrow\cr
    &&\uparrow&&&\Uparrow\cr
    &&\circ&\tto&\circ&\Rightarrow\cr
    \end{matrix} 
$$

\noindent Modules $V(p; m+n\varepsilon; dir; \mu)$ are determined by 
three nonnegative integers: the superdimension $m+n\varepsilon$ of the 
subsuperspace corresponding to each node and the total number of 
nodes, $p$; in either of the two slanted lines the graph ``lies" 
on and the directions $dir$ of the first arrow (since we only consider 
graded modules, the odd operators $a$ and $b$ must change parity; 
hence, the direction of the last noncurvy arrow must coincide with 
that of the first one), and the parameter $\mu$ corresponds to the 
eigenvalue of the Jordan cell represented by the curvy arrow (not 
depicted).

The following result is an easy corollary of \cite{GP}, pp.  59--60 
(compare with \cite{ZN}, where the description of type II modules 
contains an omission).

\begin{Lemma} Let $V$ be a reduced indecomposable $\Lambda (2)$-module.  Up 
to the change of parity, $V$ is one of the above modules 
$V(p+q\varepsilon; dir)$ of type I or modules $V(p; 
m+n\varepsilon; dir; \mu)$ of type II.
\end{Lemma}

\ssec{1.4.  Indecomposable representations of $Cl_\hbar(2)$ for 
$\hbar\neq 0$} Let $X_\pm$ be the generators of $\fsl(1|1)$.  Denote 
by $V^\hbar(n)$ for $n>0$ the $\fsl(1|1)$-module induced from the 
$n$-dimensional representation $\rho_\hbar$ of the uppertriangular 
subalgebra $\fb=Span (E, X_+)$:
$$
\rho_\hbar(E)=J_n(\hbar) \; \text{ for the Jordan cell }\; 
J_n(\hbar),\quad \rho_\hbar(X_+)=0.
$$
Clearly, $V^\hbar(1)=V^\hbar$ and on $V^\hbar(n)$ we have:
$$
\rho_\hbar(X_-)=\pmatrix 0&0\\ 1_n&0\endpmatrix,\quad 
\rho_\hbar(E)=\pmatrix J_n(\hbar)&0\\ 0& J_n(\hbar)\endpmatrix,\quad 
\rho_\hbar(X_+)=\pmatrix 0&J_n(\hbar)\\ 0& 0\endpmatrix.
$$

\begin{Lemma} Up to the change of parity, all the indecomposable 
representations of $Cl_\hbar(2)$ for $\hbar\neq 0$ are realized in the
modules $V^\hbar(n)$.
\end{Lemma}

\begin{proof} A representation from Lemma is, clearly, indecomposable.  
Such representatins realize the only way to glue two copies of 
$V^\hbar$: indeed, $\dim H^1(\fsl(1|1); \End (V^\hbar))=1$.  Since 
$\dim H^1(\fsl(1|1); V^\hbar\otimes V^{\tilde{\hbar}})=0$ for 
$\hbar\neq\tilde{\hbar}$, it is impossible to glue modules $ V^\hbar$ 
for different $\hbar$'s.  Similar argument applies to $V^\hbar(n)$.  
\end{proof}

Now, let us collect our results.

\ssbegin{1.5}{Theorem} Indecomposable finite dimensional 
representations of $\fsl(1|1)$ are realized, up o the change of 
parity, in the modules

{\em 1)} $V^\hbar(n)$ of dimension $n+n\varepsilon$;

{\em 2)} $V(p+q\varepsilon; dir)$ of type I; its dimension is equal 
to $p+q\varepsilon$ with $q=p$ or $p\pm 1$;

{\em 3)} $V(p; m+n\varepsilon; dir; \mu)$ of type II; its dimension is 
equal to $p(m+n)(1+\varepsilon)$;

{\em 4)} $\fgl(1|1)$; this is a free $\Lambda (2)$-module.
\end{Theorem}

\section*{\protect\S 2.  Irreducible representations of 
$\fg\not\simeq\fsl(1|1)$}
 
The Lie superalgebras $\fg$ of series $\fsl$ and $\fosp(2|2n)$ have a 
$\Zee$-grading of the form $\fg=\fg_{-1}\oplus\fg_{0}\oplus\fg_{1}$.  
Hence, (\cite{BL}), any irreducible $\fg$-module $L^h$ with 
highest weight $h$ and even highest vector is the quotient of the 
induced module $I(V^h)=U(\fg)\otimes_{U(\fg_{0}\oplus\fg_{1})}V^h$ for 
an irreducible $\fg_{0}$-module $V^h$ with highest weight $h$.

V.~Kac gave a description \cite{K} of conditions under which the 
module $I(V^h)$ is irreducible, if $\dim V^h<\infty$.  A weight $h$ is 
{\it typical} if $(h+\rho, \varphi_i)\neq 0$ for the Killing form 
$(\cdot, \cdot)$ on $\fg$ and every root $\varphi_i$ of the odd space 
$\fg_{-1}$.  (Here $\rho$ is a half sum of positive even roots and 
negative odd ones.)

For atypical weights the modules $I(h)$ constitute infinite in both 
directions acyclic complexes, first described in\cite{BL} and 
\cite{L1}.  Recall this result.

\ssec{2.1.  The complex of integral and differential forms and its 
generalizations} For $\fsl(1|2)\cong\fosp(2|2)\cong\fvect(0|2)$ the 
interpretation of the modules in the realization by vector fields is 
the most graphic one.  Every irreducible finite dimensional 
$\fvect(0|2)$-module is of the following form (up to $\Pi$).  Let $(a, 
b)$ be the weights with respect to $\xi_1\partial_1$ and 
$\xi_2\partial_2$, respectively.  Then $I(V^{(a, b)})$ is irreducible 
if $(a, b)\neq (0, -n)$ or $(n+1, 1)$.  In the latter case, set 
$\Sigma_{-n}=I(V^{(0, -n)})$ and $\Omega^{n}=I(V^{(n+1, 1)})$.  These 
spaces are called the superspaces of integral $(-n)$-forms and 
differential $n$-forms, respectively.  They constitute an acyclic 
complex
$$
\dots\stackrel{d}{\longrightarrow}\Sigma_{-n}\stackrel{d} 
{\longrightarrow}\dots \stackrel{d}{\longrightarrow}
\Sigma_{-1}\stackrel{d}{\longrightarrow} 
\Sigma_{0}\stackrel{\int}{\longrightarrow}\Omega^{0} 
\stackrel{d}{\longrightarrow} \dots\stackrel{d}{\longrightarrow} 
\Omega^{n}\stackrel{d}{\longrightarrow}\dots
$$
Denote by $i(k)$ the kernel of the outgoing arrow at the $k$-th place 
of the above complex.  This is an irreducible module and two modules 
with neighboring numbers are glued in one indecomposable one, more 
exactly,
$$
\Sigma_{-n}\simeq i(-n+1)\longrightarrow \Pi(i(-n))\;  \text{ and }\; 
\Omega^n\simeq\Sigma_{-n}^*\simeq i(n)\longrightarrow \Pi(i(n+1)).
$$
These arrows can be organized in graphs similar to the cases 
$\fg=\fsl(1|1)$.  Let us number the nodes in graphs corresponding to 
modules of type I downwards, starting with an integer $k$.  Let the 
$n$-th node denote the module $i(n)$ or $\Pi(i(n))$, and denote the 
corresponding module by $V(p+q\varepsilon, dir; k)$.  The dimension of 
the module obtained is equal to $2p(k+p-1)(1+\eps)$ if $p=q$ and to 
$N+(N-1)\eps$, where $N=k(2p+1)+2p^2$, for $p=q-1$, to $N+(N+1)\eps$, 
where $N=k(2p+1)+2(p^2-p+1)$, for $p=q+1$ The following module $Sq(k)$ 
of dimension $(4k-2)(1+\varepsilon)$ is an analog of the 
representation of $\fsl(1|1)$ in $\fgl(1|1)$, i.e., the free $\Lambda 
(2)$-module of rank 1:
$$
\begin{matrix}
    i(k)&\tto&\Pi(i(k+1))\cr
    \downarrow&&\downarrow\cr
    \Pi(i(k-1))&\tto&i(k)\cr
\end{matrix}
$$

\begin{Lemma}{\em(\cite{Sh1})} For $\fg\not\simeq \fsl(1|1)$ we 
have $\Ext^1(I(\varphi), I(\chi))=0$;
$$
\Ext^1(i(\varphi), i(\chi))=\cases \Cee&\text{if, up to $\Pi$,  $(i(\varphi),
i(\chi))=(I(\varphi), I(\chi))$}\\
0&\text{otherwise.}\endcases
$$ 
\end{Lemma}

This lemma shows that indecomposable modules over $\fg\not\simeq
\fsl(1|1)$ are simpler than those over $\fsl(1|1)$.

\section*{\protect\S 3. Supervarieties of representations}

So far we have described points of the {\it variety} that parametrizes 
$\fg$-modules.

To consider the supervariety is not difficult: its {\it points} are 
the same as those of the underlying variety.  The tangent space to the 
supervariety of parameters at its point corresponding to an 
indecomposiable module $V$ is, as follows from the cohomology theory 
(cf., e.g., \cite{F}), isomorphic to $H^1(\fg; \End (V))$.

\begin{Theorem}{\em (\cite{Sh1})} $H^1(\fg; \End (V))_{\od}=0$ for 
$\fg\not\cong\fsl(1|1)$ and any indecomposiable module $V$ or if 
$\fg\cong\fsl(1|1)$ and $V\cong V_\hbar(n)$.
\end{Theorem}

Proof is straightforward, with the help of the Cazimir element. \qed

For $V\not\cong V_\hbar(n)$ there are odd parameters.  Indeed, 
consider the simplest cases: $V\cong {\bf 1}$ or $\Pi({\bf 1})$.  
In either case, $\End (V)\cong{\bf 1}$ and, since 
$H^1(\fsl(1|1))\cong\Cee[\xi, \eta]/(\xi\eta)$, where $\xi$ and $\eta$ 
are odd 1-cocycles, see \cite{FL}, $\dim H^1(\fsl(1|1))_{\od}=2$.
There are no obstructions to globalization of these deformations.

Computation of $H^1(\fg; \End (V))$ for 
modules more complicated than the trivial one seems at the moment to 
be too difficult to handle by bare hands or even bare computers.

\section*{\protect\S 4. On indecomposable representations of vectorial 
Lie superalgebras}

Consider any of the simple finite dimensional vectorial Lie 
superalgebras $\cL$ except $\widetilde \fsvect(0|n)$ or $\fspe(n)$ in 
their standard $\Zee$-grading: 
$\cL=\mathop{\oplus}\limits_{-1}^kL_{i}$.  Set
$$
L_{\geq}=\mathop{\oplus}\limits_{i\geq 0}L_{i};\; 
L_{>}=\mathop{\oplus}\limits_{i>0}L_{i}; \; 
L_{-}=\mathop{\oplus}\limits_{i<0}L_{i};\;  
L_{\leq}=\mathop{\oplus}\limits_{i\leq 0}L_{i}.
$$
Let $V$ be an irreducible $L_{0}$-module.  We will consider two types 
of $\cL$-modules:

1) let $L_{>}V=0$ and set $I(V)=U(\cL)\otimes_{U(L_{\geq}}V$. 

2) let $L_{<}V=0$ and set $\cI(V)=U(\cL)\otimes_{U(L_{\leq}}V$. 

The Lie superalgebras with Cartan matrix and with a compatible 
grading are of the form $\cL=L_{-1}\oplus L_{0}\oplus L_{1}$, 
where $L_{1}\simeq L_{-1}^*$. Therefore, for them $I(V)\simeq\cI(V^*)$.

Contrarywise for the vectorial Lie superalgebras the diference between 
modules $I(V)$ and $\cI(V)$ is crucial.  As is easy to see with the help 
of the Poincar\'e-Birkhoff-Witt theorem, $I(V)\simeq \Lambda (L_{-})\otimes 
V$, as spaces.  For vectorial Lie superalgebras distinct from 
$\fsvect(0|2)$, $\fsvect(0|3)\simeq \fspe(3)$ and $\fs\fh(0|4)$ the 
space $\cI(V)$ is of infinite dimension.  

For $\fsvect(0|3)$ and $\fspe(n)$, though $\cI(V)$ is of finite
dimension, still $I(V)\not\simeq\cI(V^*)$.  However, thanks to the
existence of complete description of typical $\fspe(n)$-modules, see
\cite{L3} one may hope for a complete description of indecomposable
$\fspe(n)$-modules, at least, for $n=3$.  The limitations for this
hope are set by Theorem 4.2.

\ssbegin{4.1}{Theorem} {\em (Shapovalov, 1985)} Typical irreducible 
representaitons of $\fvect(0|2)$ and $\fs\fh(0|4)$ are direct 
summands.  \end{Theorem}

\ssbegin{4.2}{Theorem} {\em (Shapovalov, 1985)} Indecomposable 
representaitons of $\fvect(0|n)$ for $n>2$, $\fs\fpe(n)$, 
$\fsvect(0|n)$, $\widetilde\fsvect(0|n)$ for $n>3$ and $\fs\fh(0|n)$ for 
$n>4$ may include any irreducible module as the composition factor.
\end{Theorem}

No proof of these theorems was ever published. Recent result of 
Shomron (given with proof) \cite{Sh} covers the $\fvect(0|n)$ case and 
indicates how to tackle the other cases.

\section*{Appendix}

\ssec{A.1.  Linear algebra in superspaces.  Generalities} A {\it 
superspace} is a $\Zee /2$-graded space; for a superspace 
$V=V_{\ev}\oplus V_{\od }$ denote by $\Pi (V)$ another copy of the 
same superspace: with the shifted parity, i.e., $(\Pi(V))_{\bar i}= 
V_{\bar i+\od }$.  The {\it superdimension} of $V$ is $\dim V=p+q\eps 
$, where $\eps ^2=1$ and $p=\dim V_{\ev}$, $q=\dim V_{\od }$.  
(Usually $\dim V$ is expressed as a pair $(p, q)$ or $p|q$; this 
obscures the fact that $\dim V\otimes W=\dim V\cdot \dim W$ which 
becomes manifest with the use of $\eps$.)

A superspace structure in $V$ induces the superspace structure in the 
space $\End (V)$.  A {\it basis of a superspace} is always a basis 
consisting of {\it homogeneous} vectors; let $\Par=(p_1, \dots, 
p_{\dim V})$ be an ordered collection of their parities.  We call 
$\Par$ the {\it format} of the basis of $V$.  A square {\it 
supermatrix} of format (size) $\Par$ is a $\dim V\times \dim V$ matrix 
whose $i$th row and $i$th column are of the same parity $p_i$.  The 
matrix unit $E_{ij}$ is supposed to be of parity $p_i+p_j$ and the 
bracket of supermatrices (of the same format) is defined via Sign 
Rule: 

{\it if something of parity $p$ moves past something of parity 
$q$ the sign $(-1)^{pq}$ accrues; the formulas defined on homogeneous 
elements are extended to arbitrary ones via linearity}.  

Examples of application of Sign Rule: 
setting $[X, Y]=XY-(-1)^{p(X)p(Y)}YX$ we get the notion of the 
supercommutator and the ensuing notions of the supercommutative 
superalgebra and the Lie superalgebra (that in addition to 
superskew-commutativity satisfies the super Jacobi identity, i.e., the 
Jacobi identity amended with the Sign Rule). The derivation of a 
superalgebra $A$ is a linear map $D: A\tto A$ such that satisfies the 
Leibniz rule (and Sign rule)
$$
D(ab)=D(a)b+(-1)^{p(D)p(a)}aD(b). 
$$
In particular, let $A=\Kee[x]$ be the free supercommutative polynomial 
superalgebra in $x=(x_{1}, \dots , x_{n})$, where the superstructure 
is determined by the parities of the indeterminates: $p(x_{i})=p_{i}$.  
Partial derivatives are defined (with the help of Leibniz and Sign 
Rules) by the formulas $\pderf{x_{i}}{x_{j}}=\delta_{i,j}$.

Clearly, the collection $\fder A$ of all superdifferentiations of $A$ 
is a Lie superalgebra $\fvect(m|n)$ on $m$ even and $n$ odd 
indeterminates whose elements are of the form
$$
D=\sum f_i(x)\pder{x_{i}}.
$$
The {\it divergence} of such $D$ is defined to be
$$
\Div D=\sum (-1)^{p(f_{i})}\pderf{f_i}{x_{i}}.
$$
Set $\fsvect(m|n)=\{D\in \fvect(m|n)\mid \Div D=0\}$ and $\widetilde 
\fsvect(m|n)=\{D\in \fvect(m|n)\mid \Div (1+t\xi_{1}\dots 
\xi_{n})D=0\}$, where $p(t)\cong n\pmod 2$.

Define the superalgebra of differential forms as the supercommutative 
superalgebra of polynomials in the $x_{i}$ and $dx_{i}$ with 
$p(dx_{i})=p(x_{i})+1$.  Define the exterior differential $d$ by the 
formula $d(f)=\sum dx_i\pderf{f_i}{x_{i}}$ and extend $d$ to forms of 
higher degrees (in $dx$) via Leibniz and Sign Rules.  Define the Lie 
derivative $L_{D}$ along $D\in \fvect(m|n)$ by the formula 
$L_{D}d(f)=(-1)^{p(D)}d(D(f))$ extended to higher forms via Leibniz 
and Sign Rules.

Set $\omega =\sum (dx_{i})^2$ and $\fh(0|n)=\{D\in \fvect(0|n)\mid 
L_{D}(\omega)=0\}$.

The {\it general linear} 
Lie superalgebra of all supermatrices of size $\Par$ is denoted by 
$\fgl(\Par)$; usually, $\fgl(\ev, \dots, \ev, \od, \dots, \od)$ is 
abbreviated to $\fgl(\dim V_{\bar 0}|\dim V_{\bar 1})$.  Usually, 
$\Par$ is of the form $(\ev , \dots, \ev , \od , \dots, \od)$.  Such a 
format is called {\it standard}.  In this paper we can do without 
nonstandard formats but they are vital in the study of systems of 
simple roots that the reader might be interested in.  Any matrix from 
$\fgl(\Par)$ can be expressed as the sum of its even and odd parts; in 
the standard format this is the block expression:
$$
\begin{pmatrix}A&B\\ C&D\end{pmatrix}=\begin{pmatrix}A&0\\
0&D\end{pmatrix}+\begin{pmatrix}0&B\\ C&0\end{pmatrix},\quad 
p\left(\begin{pmatrix}A&0\\
0&D\end{pmatrix}\right)=\ev, \; p\left(\begin{pmatrix}0&B\\
C&0\end{pmatrix}\right)=\od.
$$

The {\it supertrace} is the map $\fgl (\Par)\longrightarrow \Cee$, 
$(A_{ij})\mapsto \sum (-1)^{p_{i}}A_{ii}$.  Since $\str [x, y]=0$, the 
space of supertraceless matrices constitutes the {\it special linear} 
Lie subsuperalgebra $\fsl(\Par)$.

{\bf Lie superalgebras that preserve bilinear forms: two types}.  To 
the linear map $F$ of superspaces there corresponds the dual map $F^*$ 
between the dual superspaces; if $A$ is the supermatrix corresponding 
to $F$ in a basis of the format $\Par$, then to $F^*$ the {\it 
supertransposed} matrix $A^{st}$ corresponds:
$$
(A^{st})_{ij}=(-1)^{(p_{i}+p_{j})(p_{i}+p(A))}A_{ji}.
$$

The supermatrices $X\in\fgl(\Par)$ such that 
$$
X^{st}B+(-1)^{p(X)p(B)}BX=0\quad \text{for an homogeneous matrix 
$B\in\fgl(\Par)$}
$$
constitute the Lie superalgebra $\faut (B)$ that preserves the 
bilinear form on $V$ with matrix $B$.  

Recall that the {\it supersymmetry} of the homogeneous form $\omega $ 
means that its matrix $B$ satisfies the condition $B^{u}=B$, where 
$B^{u}=
\begin{pmatrix} 
R^{t} & (-1)^{p(B)}T^{t} \\ (-1)^{p(B)}S^{t} & -U^{t}\end{pmatrix}$ 
for the matrix $B=\begin{pmatrix} R& S \\ T & U\end{pmatrix}$. 
Similarly, {\it skew-supersymmetry} of $B$ means that $B^{u}=-B$. 
Thus, we see that the {\it upsetting} of bilinear forms $u:\Bil (V, 
W)\tto\Bil(W, V)$, which for the {\it spaces} $V=W$ is expressed on 
matrices in terms of the transposition, is a new operation.

Most popular canonical forms of the nondegenerate supersymmetric form 
are the ones whose supermatrices in the standard format are the 
following canonical ones, $B_{ev}$ or $B'_{ev}$:
$$
B_{ev}(m|2n)= \begin{pmatrix} 
1_m&0\\
0&J_{2n}
\end{pmatrix},\quad \text{where 
$J_{2n}=\begin{pmatrix}0&1_n\\-1_n&0\end{pmatrix}$,}
$$
or
$$
B'_{ev}(m|2n)= \begin{pmatrix} 
\antidiag (1, \dots , 1)&0\\
0&J_{2n}
\end{pmatrix}. 
$$
The usual notation for $\faut (B_{ev}(m|2n))$ is $\fosp(m|2n)$ or 
$\fosp^{sy}(m|2n)$.  (Observe that the passage from $V$ to $\Pi (V)$ 
sends the supersymmetric forms to superskew-symmetric ones, preserved 
by the \lq\lq symplectico-orthogonal" Lie superalgebra $\fsp'\fo 
(2n|m)$ or $\fosp^{sk}(m|2n)$ which is isomorphic to 
$\fosp^{sy}(m|2n)$ but has a different matrix realization.  We never 
use notation $\fsp'\fo (2n|m)$ in order not to confuse with the 
special Poisson superalgebra.)

In the standard format the matrix realizations of these algebras are:
$$
\begin{matrix} 
\fosp (m|2n)=\left\{\left (\begin{matrix} E&Y&X^t\\
X&A&B\\
-Y^t&C&-A^t\end{matrix} \right)\right\};\quad \fosp^{sk}(m|2n)=
\left\{\left(\begin{matrix} A&B&X\\
C&-A^t&Y^t\\
Y&-X^t&E\end{matrix} \right)\right\}, \\
\text{where}\; 
\left(\begin{matrix} A&B\\
C&-A^t\end{matrix} \right)\in \fsp(2n),\quad E\in\fo(m)\;
\text{and}\;  {}^t \; \text{is the usual transposition}.\end{matrix} 
$$

A nondegenerate supersymmetric odd bilinear form $B_{odd}(n|n)$ can be 
reduced to the canonical form whose matrix in the standard format is 
$J_{2n}$.  A canonical form of the superskew odd nondegenerate form in 
the standard format is $\Pi_{2n}=\begin{pmatrix} 
0&1_n\\1_n&0\end{pmatrix}$.  The usual notation for $\faut 
(B_{odd}(\Par))$ is $\fpe(\Par)$.  The passage from $V$ to $\Pi (V)$ 
sends the supersymmetric forms to superskew-symmetric ones and 
establishes an isomorphism $\fpe^{sy}(\Par)\cong\fpe^{sk}(\Par)$.  
This Lie superalgebra is called, as A.~Weil suggested, {\it 
periplectic}.  In the standard format these superalgebras are 
shorthanded as in the following formula, where their matrix 
realizations is also given:
$$
\begin{matrix}
\fpe ^{sy} (n)=\left\{\begin{pmatrix} A&B\\
C&-A^t\end{pmatrix}, \; \text{where}\; B=-B^t,
C=C^t\right\};\\
\fpe^{sk}(n)=\left\{\begin{pmatrix}A&B\\ C&-A^t\end{pmatrix}, \;
\text{where}\; B=B^t, C=-C^t\right\}.
\end{matrix}
$$

The {\it special periplectic} superalgebra is $\fspe(n)=\{X\in\fpe(n): \str
X=0\}$.

Observe that though the Lie superalgebras $\fosp^{sy} (m|2n)$ and 
$\fpe ^{sk} (2n|m)$, as well as $\fpe ^{sy} (n)$ and $\fpe ^{sk} (n)$, 
are isomorphic, the difference between them is sometimes crucial.

\end{document}